\documentclass[letterpaper, 10 pt, conference]{ieeeconf}  % Comment this line out
\IEEEoverridecommandlockouts                              % This command is only
\usepackage{ulem}
\usepackage[backref,colorlinks=true, linkcolor=blue, citecolor=blue, urlcolor=blue, pdfborder={0 0 0}]{hyperref}
\usepackage{amssymb,amsmath}
\usepackage{bbm}
\usepackage{graphicx,comment}
\usepackage{algorithm,algorithmic}
\usepackage{hyperref}
\usepackage{colortbl}
\usepackage{authblk}
\usepackage{subcaption}
\newcommand{\R}{\mathbb{R}}

\newcommand{\C}{\mathbb{C}}
\newcommand{\F}{\mathcal{F}}
\newcommand{\re}{\mathrm{Re}}

\usepackage{url}
\graphicspath{{images/}}
\usepackage[symbol]{footmisc}

\title{\LARGE \bf
Generalized Proximal Smoothing for Phase Retrieval
}

\author{Minh Pham$^{1}$, Penghang Yin$^{1}$, Arjun Rana$^{2}$, Stanley Osher$^{1}$ and Jianwei Miao$^{2}$% <-this % stops a space
\thanks{$^{1}$Department of Mathematics, University of California, Los Angeles, Los Angeles, CA 90095, USA, \textit{Email: (minhrose,yph,sjo)@math.ucla.edu}}%
\thanks{$^{2}$Department of Physics, University of California, Los Angeles, Los Angeles, CA 90095, USA, \textit{Email: (arana,miao)@physics.ucla.edu }}%
}

\begin{document}

\maketitle
\thispagestyle{empty}
\pagestyle{empty}

%%%%%%%%%%%%%%%%%%%%%%%%%%%%%%%%%%%%%%%%%%%%%%%%%%%%%%%%%%%%%%%%%%%%%%%%%%%%%%%%
\begin{abstract}

In this paper, we report the development of the generalized proximal smoothing (GPS) algorithm for phase retrieval of noisy data. GPS is a optimization-based algorithm, in which we relax both the Fourier magnitudes and object constraints. We relax the object constraint by introducing the generalized Moreau-Yosida regularization and heat kernel smoothing. We are able to readily handle the associated proximal mapping in the dual variable by using an infimal convolution. We also relax the magnitude constraint into a least squares fidelity term, whose proximal mapping is available. GPS alternatively iterates between the two proximal mappings in primal and dual spaces, respectively. Using both numerical simulation and experimental data, we show that GPS algorithm consistently outperforms the classical phase retrieval algorithms such as hybrid input-output (HIO) and oversampling smoothness (OSS), in terms of the convergence speed, consistency of the phase retrieval, and robustness to noise.


\end{abstract}

%%%%%%%%%%%%%%%%%%%%%%%%%%%%%%%%%%%%%%%%%%%%%%%%%%%%%%%%%%%%%%%%%%%%%%%%%%%%%%%%
\section{Introduction}

Phase retrieval has been fundamental to several disciplines, ranging from imaging, microscopy, crystallography and optics to astronomy \cite{Millane1990,miao2015beyond,shechtman2015phase,marchesini2007invited}. It aims to recover an object only from its Fourier magnitudes. Without the Fourier phases, the recovery can be achieved via iterative algorithms when the Fourier magnitudes are sampled at a frequency sufficiently finer than the Nyquist interval \cite{Miao1998}. In 1972, Gerchberg and Saxton developed an iterative algorithm for phase retrieval, utilizing the magnitude of an image and the Fourier magnitudes as constraints \cite{Gerchberg}. In 1982, Fienup generalized the Gerchberg-Saxton algorithm by developing two iterative algorithms: error reduction (ER) and hybrid input-output (HIO), which use a support and positivity as constraints in real space and the Fourier magnitudes in reciprocal space \cite{fienup1978reconstruction}. In 1998, Miao, Sayre and Chapman proposed, when the number of independently measured Fourier magnitudes is larger than the number of unknown variables associated with a sample, the phases are in principle encoded in the Fourier magnitudes and can be retrieved by iterative algorithms \cite{Miao1998}. These developments finally led to the first experimental demonstration of coherent diffractive imaging (CDI) by Miao and collaborators in 1999 \cite{Miao1999}, which has stimulated wide spread research activities in phase retrieval, CDI, and their applications in the physical and biological sciences ever since \cite{miao2015beyond,chapman2010coherent,robinson2009coherent}.

For a finite object, when its Fourier transform is sampled at a frequency finer than the Nyquist interval (i.e. oversampled), mathematically it is equivalent to padding zeros to the object in real space. In another words, when the magnitude of the Fourier transform is oversampled, the correct phases correspond to the zero-density region surrounding the object, which is known as the oversampling theorem \cite{Miao1998,miao2000possible}. The phase retrieval algorithms iterate between real and reciprocal space using zero-density region and the Fourier magnitudes  as  dual-space constraints. A support is typically defined to separate the zero-density region from the object. The positivity constraint is applied to the density inside the support.  In the ER algorithm, the no-density region outside the support and the negative density inside the support are set to zero in each iteration \cite{fienup1978reconstruction}. The HIO algorithm  relaxes the ER in the sense that  it gradually reduces  the densities that  violate the object constraint instead  of directly forcing them  to zero \cite{fienup1978reconstruction}. This relaxation often leads to good reconstructions from noise-free patterns. However, in real experiments, the diffraction intensities, which are proportional to the square of Fourier magnitudes, are corrupted  by a combination of  Gaussian and Poisson noise and missing data. In the presence of experimental  noise and missing data, phase retrieval  becomes  much more challenging, and the ER and HIO algorithms may only converge to sub-optimal  solutions.  Simply combining ER and HIO still suffers from stagnation and the iterations can get trapped  at local minima \cite{Fienup}.  To alleviate these problems, more advanced phase retrieval algorithms have been developed such as the shrink-wrap algorithm and guided HIO (gHIO) \cite{chen2007application,marchesini2007invited}. In 2010, a smoothness constraint in real space was first introduced to improve  the phase retrieval  of noisy data \cite{Raines2010}.   Later, a noise robust framework was implemented  for enhancing the performance of existing algorithms \cite{Martin2012}. Recently, Rodriguez et al. proposed to impose the smoothness constraint on the no-density region outside the support by applying Gaussian filters \cite{Rodriguez}.  The resulting oversampling smoothness  (OSS)  algorithm  successfully reduces  oscillations  in the  reconstructed  image, and is more robust  to noisy data than the existing algorithms.

Since phase retrieval can be cast as a non-convex minimization problem, many efforts have been made to study phase retrieval algorithms from the viewpoint of optimization. For example, Bauschke et al. \cite{Bauschke2} related HIO to a particular relaxation of the Douglas-Rachford algorithm \cite{Douglas} and introduced the hybrid projection reflection algorithm \cite{Bauschke1, Bauschke2}. Using similar ideas, researchers further proposed several projection algorithms such as iterated difference map \cite{Elser} and relaxed averaged alternation reflection \cite{Luke}. In \cite{Fannjiang_16}, Chen and Fannjiang analyzed a Fourier-domain Douglas-Rachford algorithm for phase retrieval. By taking noise into account, the Wirtinger Flow \cite{CandesWF} relaxes the magnitude constraint into a fidelity term that measures the misfit of Fourier data, to which Wirtinger gradient descent is applied. Other methods in this line include alternating direction methods \cite{Chang1,Chang2,ADM_wen} that have been widely used in image processing, as well as lifting approaches \cite{Chai_12} such as PhaseLift \cite{CandesPL,CandesPL2} by Cand\`es et al. and its variants \cite{Waldspurger,Yin2015}.

In this paper, we propose an optimization-based phase retrieval method, termed generalized proximal smoothing (GPS), which effectively addresses the noise in both real and Fourier spaces. Motivated by the success of OSS \cite{Rodriguez}, GPS incorporates the idea of Moreau-Yosida \cite{Moreau_65,Yosida_64} regularization with heat kernel smoothing, to relax the object constraint into an implicit regularizer. We extend the notion of infimal convolution from real domain to complex domain in the context of convex analysis \cite{rockafellar2015convex,boyd2004convex}, which enables us to handle the convex conjugate of the implicit relaxation in the dual variable. We further relax the magnitude constraint into a least squares fidelity term, for de-noising in Fourier space. To minimize the primal-dual formulation, GPS iterates back and forth between efficient proximal mappings of the two relaxed functions, respectively. Our experimental results using noisy experimental data of biological and inorganic specimens demonstrate that GPS consistently outperforms the state-of-the-art algorithms HIO and OSS in terms of both speed and robustness. We also refer readers to the recent paper \cite{binaryrelax} about training quantized neural networks, which shows another success of using Moreau-Yosida regularization to relax the hard constraint.\\

\noindent \textbf{Notations.} Let us fix some notations. For any complex-valued vectors $u, v \in \mathbb{C}^{ n}$, $\overline{u}$ is the complex conjugate of $u$, whereas $u^* := \overline{u}^\top$ is the Hermitian transpose. $\re(u)$ and $\mathrm{Im}(u)$ are the real and imaginary parts of $u$, respectively.
$$
\langle u , v\rangle := u^* v = \sum_{i=1}^n \overline{u_{i}} v_{i}
$$ 
is the Hermitian inner product of $u$ and $v$. $u \circ v$ is the element-wise product of $u$ and $v$ given by $(u\circ v)_{i} = u_{i} v_{i}$. $\| u \| := \sqrt{\langle u , u \rangle}$ denotes the $\ell_2$ norm of $u$. Given any Hermitian positive semi-definite matrix $K\in\C^{n\times n}$, we define $\| u \|_{K} := \sqrt{\langle u , K \,u \rangle}$. $\textrm{arg}(u)$ denotes the argument (or phase) of $u = \re(u)+\mathbf{i}\cdot\mathrm{Im}(u)$, which is given by
\begin{align*}
	\textrm{arg}(u) := \begin{cases}
    	\tan^{-1}\big( \frac{\mathrm{Im}(u)}{\re(u)} \big)	& \mbox{if } u \neq 0, \\
        1		& \mbox{otherwise.}
    \end{cases}
\end{align*} 
$\mathcal{I}_{\mathcal{X}}$ is the characteristic function of a closed set $\mathcal{X}\subset\C^n$ given by
    \begin{align*}
        \mathcal{I}_{\mathcal{X}}(x) = \begin{cases}
            0 & x \in \mathcal{X} \\
            \infty & \mbox{otherwise.}
        \end{cases}
    \end{align*}
$\mathrm{proj}_{\mathcal{X}}(u):= \arg\min_{v\in\mathcal{X}} \|v-u\|$ 
is the projection of $u$ onto $\mathcal{X}$, and
$\mathop{\mathrm{prox}}_{f}$ is the proximal mapping of the function $f(u)$ defined by
$$
\mathrm{prox}_{f}(u) := \arg\min_v \; \Big\{ f(v) + \frac{1}{2}\|v-u\|^2 \Big\}.
$$

\bigskip

\section{Proposed Model}
We consider the reconstruction of a 2D image $u$ defined on a discrete lattice 
$$
\Omega:=\{(i,j): 1\leq i\leq n_1, 1\leq j\leq n_2\}.
$$
For simplicity, we represent $u$ in terms of a vector in $\R^n$ by the lexicographical order with $n = n_1\times n_2$. Then $u_{i}$ represents the density of image at the $i$-th pixel. Due to oversampling, the object densities reside in a sub-domain $S\subset\Omega$ known as the support, and $u$ is supposed to be zero outside $S$. Throughout the paper, we assume that the support $S$ is centered around the domain $\Omega$. The object constraint is
\begin{align*}
\mathcal{S}:= \{u\in\R^{n}: u_{i}\geq0 \mbox{ if }  i\in S, \; u_{i}=0 \mbox{ otherwise} \}.
\end{align*} 
The Fourier magnitude data is obtained as $b = |\mathcal{F} u|$, where $\mathcal{F}:\mathbb{R}^{n_1\times n_2}\to\mathbb{C}^{n_1\times n_2}$ is the discrete Fourier transform (DFT). We denote the magnitude constraint by
$$\mathcal{T}: = \{u \in \mathbb{R}^n: |\mathcal{F}u | = b \}.$$
In the absence of noise, phase retrieval (PR) problem is simply to
$$
\mathrm{find} \; u\in \R^{n}, \quad \mbox{such that } u\in\mathcal{S} \cap \mathcal{T}.
$$ 
This amounts to the following composite minimization problem
\begin{align}\label{eq:pr_comp}
    \min_{u \in\R^{n}} \; f(u) +  g(\F u), 
\end{align}
where $f(u):= \mathcal{I}_{\mathcal{S}}(u)$ and $g(z):= \mathcal{I}_{|z|=b}(z)$ are two characteristic functions that enforce the object and Fourier magnitudes constraints. Note that $f$ is a closed and convex function while $g$ is closed but non-convex, which give the non-convex optimization problem of (\ref{eq:pr_comp}).

\subsection{A new noise-removal model}
In real experiments, the Fourier data are contaminated by experimental noise. Moreover, the densities outside the support are not exactly equal to zero either. In the noisy case, the image to be reconstructed no longer fulfills either the Fourier magnitudes or the object constraint. The ER algorithm, which alternatively projects between these two constraints, apparently does not take care of the noise. The HIO  ``relaxes'' the object constraint on densities wherever it is violated. This relaxation only helps in the noiseless case. In the presence of noise, the feasible set $\mathcal{S} \cap \mathcal{T}$ can be empty, and alternating projection methods like ER and HIO may fail to converge and keep oscillating. The OSS \cite{Rodriguez} improves the HIO by applying extra Gaussian filters to smooth the densities outside the support at different stages of the iterative process. None of them, however, seems to properly address the corruption of the Fourier magnitudes.

\medskip

Introducing the splitting variable $z = \mathcal{F} u \in\C^{n}$, we reformulate (\ref{eq:pr_comp}) as 
\begin{align}\label{eq:pr_split}
    \min_{u,z \in\C^{n}} \; f(u) +  g(z) \quad \mbox{subject to} \quad z = \mathcal{F} u.
\end{align}
Note that we need to extend $u$ to complex domain in this setting. In the presence of noise, we seek to relax the characteristic functions $f$ and $g$ that enforce hard constraints into soft constraints. To this end, we extend the definition of the Moreau-Yosida regularization \cite{Moreau_65,Yosida_64} to complex domain. Let $K$ be a Hermitian positive definite $n\times n$ matrix. The Moreau-Yosida regularization of a lower semi-continuous extended-real-valued function $h: \C^n \to (-\infty,\infty]$, associated with $K$, is defined by
$$
h_K(u) := \inf_{v\in\C^n} \; \Big\{ h(v) + \frac{1}{2}\| v - u\|_{K^{-1}}^2 \Big\}.
$$
We see that $h_K$ converges pointwise to $h$ as $\|K\|\to 0^{+}$. In the special case where $K = t I$ is a multiple of identity matrix with $t>0$, $h_K$ reduces to
$$
h_t(u) := \inf_{v\in\C^n} \; \Big\{ h(v) + \frac{1}{2t} \|v - u\|^2 \Big\}.
$$
For any characteristic function $h = \mathcal{I}_\mathcal{X}$ of a closed set $\mathcal{X}\subset \C^n$, 
$$
h_t(u) = \frac{1}{2t} \inf_{v\in\mathcal{X}} \;\|v - u\|^2 = \frac{1}{2t}\|u - \mathrm{proj}_\mathcal{X}(u)\|^2
$$
is $\frac{1}{2t}$ of the squared $\ell_2$ distance from $u$ to the set $\mathcal{X}$. Similar idea of relaxing a characteristic function into a distance function has been successfully applied to the quantization problem of deep neural networks in \cite{binaryrelax}.
Taking $\mathcal{X} = \{z\in\R^{n}: |z| = b\}$ to be the magnitude constraint set and $\sigma>0$, we first relax $g = \mathcal{I}_{|z|=b}$ in (\ref{eq:pr_split}) into
\begin{align*}
g_\sigma(z) = \frac{1}{2\sigma} \inf_{|v|=b} \; \|v-z\|^2.
\end{align*}
Since $\mathrm{proj}_{|z|=b}(z) = b\circ \exp(\mathbf{i}\cdot\arg(z))$
is the projection of $z$ onto the set $\{z\in\R^{n}: |z| = b\}$, a simple calculation shows that
\begin{align}\label{eq:g_relax}
g_\sigma(z) = & \, \frac{1}{2\sigma}\|b\circ \exp(\mathbf{i}\cdot\arg(z)) - z\|^2 \notag \\
= & \, \frac{1}{2\sigma} \|(b - |z|)\circ\exp(\mathbf{i}\cdot\arg(z))\|^2 \notag \\
= & \, \frac{1}{2\sigma}\||z| - b\|^2
\end{align}
is a least squares fidelity, which measures the difference between the observed magnitudes and fitted ones. This fidelity term has been considered in the literature by assuming the measurements being corrupted by i.i.d. Gaussian noise; see \cite{Chang1} for example. In practice, we observe that it works well even with a combination of Gaussian and Poisson noises.

\medskip

Following this line, we further relax $f = \mathcal{I}_{\mathcal{S}}$ into
\begin{align}\label{eq:f_relax}
f_{G}(u) = & \, \inf_v \; \Big\{ f(v) + \frac{1}{2} \| v - u\|_{G^{-1}}^2  \Big\} \notag\\
 = & \, 
\inf_{v\in\mathcal{S}} \; \frac{1}{2} \| v - u\|_{G^{-1}}^2
\end{align}
for some Hermitian positive definite matrix $G$. The choice of $G$ here is tricky, and will be discussed later in section \ref{sec:gps}. The relaxation of both constraints thus leads to the proposed noise-removal model
\begin{align}\label{eq:gps}
    \min_{u,z\in \mathbb{C}^n} \; f_G(u) + g_{\sigma}(z) \quad \mbox{subject to} \quad z = \mathcal{F} u.
\end{align}
For a non-diagonal matrix $G$, the associated Moreau-Yosida regularization $f_G$ in (\ref{eq:f_relax}) does not enjoy an explicit expression in general. This poses a challenge to the direct minimization of (\ref{eq:gps}) using solvers such as alternating direction method of multipliers (ADMM) \cite{glowinski1975approximation,boyd2011distributed,yan_16}.

\subsection{Generalized Legendre-Fenchel transformation}
We can express any function $h: \C^n\to \R$ as a function $\tilde{h}$ defined on $\R^{2n}$ in the following way
$$
h(u)  = \tilde{h}\big(\mathrm{Re}(u), \mathrm{Im}(u)\big) = \tilde{h}(\tilde{u}),
$$ 
where $\tilde{u} = \begin{bmatrix}
    \mathrm{Re}(u) \\
    \mathrm{Im}(u)
\end{bmatrix}\in\R^{2n}$ and $\tilde{h}:\R^{2n}\to\R$. We define that $h$ is convex, if $\tilde{h}$ is convex on $\R^{2n}$. Note that for any $u,y\in \C^n$,
$$
\re \langle u, y \rangle = \langle\re(u), \re(y) \rangle + \langle \mathrm{Im}(u), \mathrm{Im}(y)\rangle =  \langle \tilde{u}, \tilde{y} \rangle.
$$
We propose to generalize the Legendre-Fenchel transformation (a.k.a. convex conjugate) \cite{rockafellar2015convex} of an extended-real-valued  convex function $h$ defined on $\C^n$ as 
\begin{align*}
	h^*(y) := \sup_{u\in\C^n} \; \Big\{ \mathrm{Re} \langle y, u \rangle - h(u) \Big\}. 
%= \sup_{\tilde{u}\in\R^{2n}} \; \Big\{ \langle \tilde{y}, \tilde{u} \rangle - \tilde{h}(\tilde{u}) \Big\}.
\end{align*}

\medskip

In fact, $f_G$ is the infimal convolution \cite{phelps_91} between the convex functions $f$ and $\frac{1}{2}\|\cdot \|_{G^{-1}}^2$ in the sense that
\begin{equation*}
    f_G(u) = \inf_{v+w = u} \; \left\{ f(v) + \frac{1}{2} \|w \|_{G^{-1}}^2 \right\}.
\end{equation*}
Similar to the real case, the infimal convolution holds the property that
\begin{align}\label{eq:fG_conj}
f_G^*(y) = f^*(y) + \Big(\frac{1}{2}\|\cdot \|_{G^{-1}}^2 \Big)^*(y),
\end{align}
where $\Big(\frac{1}{2}\|\cdot \|_{G^{-1}}^2 \Big)^*(y) = \frac{1}{2}\|y\|_G^2$ and 
\begin{equation}\label{eq:f_conj}
   f^*(y) = \begin{cases}
        0 & \mbox{if } \mathrm{Re}(y) \le 0 \textrm{ on } S, \\
        \infty & \mbox{otherwise.}
    \end{cases} 
\end{equation}
While $f_G$ takes an implicit form, its generalized convex conjugate is readily explicit. This suggests us look at the primal-dual formulation of model (\ref{eq:gps}).

\subsection{A primal-dual formulation}
With slight abuse of notation, we say $y\in\partial h(u)$ is a subgradient of $h$ at $u$, if $\tilde{y}\in\partial{\tilde{h}(\tilde{u})}$. Then the Lagrangian of (\ref{eq:gps}) reads
\begin{equation}\label{eqn1}
    \mathcal{L}(u,z;y) =  f_G(u) + g_\sigma(z) + \mathrm{Re}\langle y, \F^* z- u\rangle,
\end{equation}
where $\F^* = \F^{-1}$ is the adjoint of $\F$ or the inverse DFT. The corresponding Karush-Kuhn-Tucker (KKT) condition is
\begin{align} \label{eqn2}
    y \in \partial f_G(u), \; -\F y \in \partial g_\sigma(z). 
\end{align}
We apply the convex conjugate and rewrite (\ref{eqn2}) as
\begin{align*}
    \F^* z = u  \in \partial f_G^*(y), \; -\F y \in \partial g_\sigma (z),
\end{align*}
which is exactly the KKT condition of the following min-max saddle point problem
\begin{equation}\label{eq:pd}
    \min_z \max_y \; g_\sigma(z) - f_G^*(y) + \mathrm{Re}\langle z, \F y\rangle,
\end{equation}
with $g_\sigma(z)$ and $f_G^*(y)$ explicitly available from equations (\ref{eq:g_relax}) and (\ref{eq:fG_conj}).

\medskip

\section{Generalized Proximal Smoothing}\label{sec:gps}
We carry out the minimization of the saddle point problem (\ref{eq:pd}) by a generalized primal dual hybrid gradient (PDHG) algorithm \cite{Chambolle2010,esser2010general,goldstein2013adaptive,Connor2017}, which iterates
\begin{equation*} 
\begin{cases}
 z^{k+1} = \mathrm{prox}_{t g_{\sigma}} \big( z^k - t \mathcal{F}y^k \big) \\
    y^{k+1} = \mathrm{prox}_{s f_G^*} \big( y^k + s\mathcal{F}^{-1} ( 2z^{k+1} - z^k) \big)
\end{cases}
%    & = \mathop{\mathrm{argmin}}_{y \in \mathcal{X}^*} \; \left\{ \frac{1}{2} \langle  y, G y \rangle +  \frac{1}{2s} \| y - \big( y^k + s\mathcal{F}^{-1} (2z^{k+1} - z^k) \big) \|^2 \right\} \nonumber
\end{equation*}
for some step sizes $s,t >0$.
The update of $z^{k+1}$ calls for computing the proximal mapping of $t g_\sigma$ \cite{Chang1}, whose analytic solution is given by 
\begin{align*}
\mathrm{prox}_{tg_{\sigma}}(z)  = & \, \arg\min_{v\in\C^n} \; g_\sigma(v) + \frac{1}{2t}\|v- z\|^2 \\
= & \, \arg\min_{v\in\C^n} \; \frac{1}{2\sigma}\||v| - b \|^2 + \frac{1}{2t}\|v- z\|^2 \\
= & \, \frac{b \circ \exp(\mathbf{i} \cdot\mathrm{arg}(z)) + (\sigma/t) z }{1+ \sigma/t}, 
\end{align*}
which is essentially a linear interpolation between $z$ and its projection onto the magnitude constraint $\{z\in\C^n: |z| = b\}$.

\medskip

Moreover, we need to find the proximal mapping of $sf^*_G$ for updating $y^{k+1}$, which reduces to
\begin{align}\label{eq:prox_fG}
\mathrm{prox}_{sf_G^*}(y) = & \, \arg\min_{v\in\C^n} \; f^*_G(v) + \frac{1}{2s}\|v- y\|^2 \notag\\
= & \, \arg\min_{v\in\C^n} \; f^*(v) + \frac{1}{2}\|v\|^2_G + \frac{1}{2s}\|v- y\|^2 \notag\\
= & \, \arg\min_{v\in\mathcal{S}^*} \; \frac{1}{2}\|v\|^2_G + \frac{1}{2s}\|v- y\|^2.
\end{align}
In the third equation, $\mathcal{S}^* := \{y\in\C^n: \re(y)\leq 0 \mbox{ on } S\}$ according to (\ref{eq:f_conj}), whose projection is
\begin{align} \label{eq:proj_Sc}
	\mathrm{proj}_{\mathcal{S}^*}(y)_i = 
    \begin{cases}
    	\re(y_i)^- + \mathbf{i}\cdot \mathrm{Im}(y_i) 	& \mbox{if } i \in S, \\
        y_i				  & \mbox{otherwise.} 
    \end{cases}
\end{align}
Here $x^- := \min(0,x)$ for $x\in\R$. Problem (\ref{eq:prox_fG}) seems to have closed-form solution only when $G$ is a diagonal matrix.

\medskip

We devise two versions of GPS algorithm based on different choices of G. Here we remark that $G$ only needs to be positive semi-definite in (\ref{eq:f_relax}), as $f_G$ can take the extended value $\infty$. In this case, although $G^{-1}$ does not exist in (\ref{eq:f_relax}), since $f_G$ is convex and lower semi-continuous, and the strong duality $f_G^{**} = f_G$ holds here, we can re-define $f_G$ via the biconjugate as
\begin{align}\label{eq:fG}
f_G(u) = & \, f_G^{**}(u) = \big(f_G^*(y)\big)^* \notag \\
= & \, \sup_{y\in\C^n} \; \Big\{ \mathrm{Re} \langle u, y \rangle - f_G^*(y)\Big\} \notag \\
= & \,  \sup_{y\in\C^n} \; \Big\{ \mathrm{Re} \langle u, y \rangle - f^*(y) - \frac{1}{2}\|y\|_G^2 \Big\} \notag \\
= & \, \sup_{y\in\mathcal{S^*}} \; \Big\{ \mathrm{Re} \langle u, y \rangle - \frac{1}{2}\|y\|_G^2 \Big\}. 
\end{align}
The remaining challenge is to solve the proximal problem (\ref{eq:prox_fG}).\\

\subsection{Real space smoothing}
One choice of $G$ is $G = \gamma \, D^\top D$. Here $D$ is the discrete gradient operator, and then $D^\top D$ is the negative of discrete Laplacian. In this case,
\begin{equation*}
    f_G^*(y) = f^*(y) + \frac{\gamma}{2} \|Dy\|^2,
\end{equation*}
which we shall refer as the real space smoothing. Since $G$ is not diagonal, the closed-form solution to (\ref{eq:prox_fG}) is not available. For small $\gamma$, we approximate the solution by
\begin{align}\label{eq:prox_fG_approx}
    \mathrm{prox}_{sf^*_G}(y) \approx ( I + s \gamma D^\top D )^{-1} \mathrm{proj}_{\mathcal{S}^*}(y). 
\end{align}
The projection is followed by the matrix inversion to ensure the smoothness of the reconstructed image after each iteration. In fact, the real space smoothing is related to the diffusion process. Consider the heat equation with an initial value condition on $\R^n\times[0,+ \infty)$:
\begin{align*}
    u_t(x,t)  = \Delta u(x,t), \; u(x,0)  =  h(x).  
\end{align*}
A numerical approach to the above problem is the backward Euler scheme:
\begin{equation*}
    u^{k+1} = \big( I + dt \, D^\top D \big)^{-1} u^k,
\end{equation*}
where $dt$ is the step size for time discretization. On the other hand, the exact solution of the heat equation is given by a Gaussian convolution of $h$
\begin{align*}
    u(x,t) = & \, \mathcal{G}_t * h = \int\mathcal{G}_t (v-x) h(v) dv, 
\end{align*}
where $\mathcal{G}_t(x) := \frac{1}{(4 \pi t)^{n/2}}\exp\big(- \frac{\|x\|^2}{4t}\big)$ is a heat kernel. This observation leads to a fast approximated implementation of (\ref{eq:prox_fG_approx}) when $\gamma$ is small:
\begin{align*}
    y^{k+1} = \mathcal{G}_{\gamma} *  \mathrm{proj}_{\mathcal{S^*}} \big( y^k + s\mathcal{F}^{-1} (2z^{k+1} - z^k) \big)
,\end{align*}
where the convolution can be done via the efficient DFT. In the context of physics, this is known as the low-pass filtering.
%$\mathcal{G}_{\gamma} * y = \re \Big( \mathcal{F}^* \mbox{ifftshift} \big[ \mbox{fftshift}  \big( \mathcal{F} y \big)  \circ \exp(-2\gamma r^2) \big] \Big)$. 
\medskip

Algorithm \ref{alg1}: GPS-R features low-pass filters for smoothing. Here we abuse notation $\gamma$ to imply the filter. Inspired by OSS \cite{Rodriguez}, we choose an increasing sequence of spatial frequency filters $\{\gamma_l\}$ (a sequence of finer filters). In our experiments, we do 1000 iterations with totally 10 stages. Each stage contains 100 iterations, in which we stick with the same filter frequency. We monitor the R-factor (relative error) during the iterative process, which is defined as
\begin{equation*}
    R_F(u^k) = \frac{\sum_i \big||\mathcal{F}u^k|_i - b_i \big|}{\sum_i b_i}.
\end{equation*}
The reconstruction with minimal $R_F$ at each state is fed into the next stage. By applying the smoothing on the entire domain, GPS-R can remove noise in real space and obtain the spatial resolution with fine features.

\begin{algorithm}
\caption{GPS-R: GPS with smoothing in real space.}
\label{alg1}
\textbf{Input}: measurements $b$, regularization parameters $\{\gamma_l\}_{l=1}^{10}$, step sizes $s$, $t >0$\\
\textbf{Initialize}: $z^0$, $y^0$.\\ 
 $\qquad R_F^{best}=1$, $z_{best} = z^0$ 
\begin{algorithmic}
\FOR {$l = 1,\dots,10$}
\STATE $y^0 = y_{best}$, $z^0 = z_{best}$
\FOR {$k = 1,\dots, 100$}
\STATE $z^{k+1} = \textrm{prox}_{tg} \big( z^k - t\mathcal{F}y^k \big)$
\STATE $y^{k+1/2} = \textrm{proj}_{\mathcal{S}^*} \big( y^k + s \mathcal{F}^{-1} \big(2z^{k+1} - z^k \big) \big)$
\STATE $y^{k+1} = \mathcal{G}_{\gamma_l} * y^{k+1/2}$
\STATE \textbf{if} $R_F^k < R_F^{best}$, \textbf{then} $R_F^{best} = R_F^k$, $z_{best} = z^k$, \textbf{end if}
\ENDFOR
\ENDFOR
\end{algorithmic}
\textbf{Output}: $\F^{-1}z_{best}$
\end{algorithm}

\subsection{Fourier space smoothing}
Another simple choice is the diagonal matrix 
$$
G = \gamma \, \textrm{Diag}(r\circ r),
$$
where $r\in\R^n$ and $r_i$ is the distance in the original 2D lattice between the $i$-th pixel and the center of image. Note that $G$ is not invertible since $r_i=0$ for the pixel at the center. By (\ref{eq:fG}),
\begin{align*}
    f_G(u) =\begin{cases}
       % \displaystyle    
        \sum_{i} \frac{ \big| u - \mathrm{proj}_{\mathcal{S}} (u) \big|^2_{i}   }{2\gamma\, r^2_{i} }& \mbox{if }  u_{i} \ge 0 \mbox{ at the center,} \\
        \infty  & \mbox{otherwise.} \\
    \end{cases}
\end{align*}
So $f_G(u)$ is a weighted sum of squares penalty on $u$. The weight is inversely proportional to the squared radius, which is infinity for density in the center. The further the density off the center, the smaller the penalty for the object constraint being violated.

\medskip

By the Parseval's identity, for square-integrable function $u$, we have
\begin{equation*}
    \int \Big|\frac{d}{d\xi} \hat{u}(\xi) \Big|^2 d\xi = \int |x \, u(x)|^2 dx,
\end{equation*}
where $\hat{u}$ is the Fourier transform of $u$. In the discrete setting, this amounts to
\begin{equation*}
    \|D \, \mathcal{F}u\|^2 = \| r \circ u\|^2,
\end{equation*}
Therefore, by (\ref{eq:fG_conj}),
\begin{equation*}
    f_G^*(y) = f^*(y) +  \frac{\gamma}{2} \|r \circ y \|^2 = f^*(y) + \frac{\gamma}{2} \|D \, \mathcal{F}y\|^2.
\end{equation*}
This means that we are smoothing $f^*(y)$ by regularizing with the $\ell_2$ gradient of Fourier coefficients of $y$. We thus refer it as the Fourier space smoothing.

\medskip

Since $G$ is diagonal, (\ref{eq:prox_fG}) has the closed-form solution
\begin{align*}
\mathrm{prox}_{sf^*_G}(y) = & \, \frac{1}{1 + s \gamma r^2} \circ \mathrm{proj}_{\mathcal{S}^*}(y) \\
\approx & \, \exp(-s\gamma r^2) \circ \mathrm{proj}_{\mathcal{S}^*}(y).
\end{align*}
The solution can be also approximated by a direct multiplication with the Gaussian filter when $\gamma$ is small. Hence, we update $y^{k+1}$ as
$$
y^{k+1} = \exp(-s\gamma r^2) \circ \mathrm{proj}_{\mathcal{S}^*}\big( y^k + s\mathcal{F}^{-1} (2z^{k+1} - z^k) \big).
$$
% Similar to GPS-R, instead of implementing the above update, we apply the heat kernel component-wise
% $$
% y^{k+1} =  \exp(-s\gamma r^2) \circ \mathrm{proj}_{\mathcal{S}^*}\big( y^k + s\mathcal{F}^{-1} (2z^{k+1} - z^k) \big).
% $$
GPS with smoothing in Fourier space (GPS-F) is summarized in Algorithm \ref{alg2}. 

\begin{algorithm}
\caption{GPS-F: GPS with smoothing in Fourier space.}
\label{alg2}
\textbf{Input}: measurements $b$, regularization parameters $\{\gamma_l\}_{l=1}^{10}$, step sizes $s, t>0$.\\
\textbf{Initialize}: $z^0$, $y^0$.\\ 
 $\qquad R_F^{best}=1$, $z_{best} = z^0$ 
\begin{algorithmic}
\FOR {$l = 1,\dots,10$}
\STATE $y^0 = y_{best}$, $z^0 = z_{best}$
\FOR {$k = 1,\dots, 100$}
\STATE $z^{k+1} = \textrm{prox}_{tg_\sigma} \big( z^k - t\mathcal{F}y^k \big)$
\STATE $y^{k+1/2} = \textrm{proj}_{\mathcal{S}^*} \big( y^k + s \mathcal{F}^{-1} (2z^{k+1} - z^k) \big)$
\STATE $\displaystyle y^{k+1} =  \exp(-s\gamma_l r^2) \circ y^{k+1/2}$
\STATE \textbf{if} $R_F^k < R_F^{best}$, \textbf{then} $R_F^{best} = R_F^k$, $z_{best} = z^k$, \textbf{end if}
\ENDFOR
\ENDFOR
\end{algorithmic}
\textbf{Output}: $\F^{-1}z_{best}$
\end{algorithm}

\subsection{Real-Fourier space smoothing}
Combining both Fourier and real space smoothing is an option. In each iteration, one can first apply the low-pass filter and then the Gaussian kernel in a heuristic way (GPS-RF).

\subsection{Incomplete measurements}
In practice,  not all diffraction intensities can be experimentally measured. For example, to prevent a detector from being damaged by an intense X-ray beam, either a beamstop has to be placed in front of the detector to block the direct beam or a hole has to be made at the center of the detector, resulting in missing data at the center \cite{Miao2005}. Furthermore, missing data may also be present in the form of gaps between detector panels.  For incomplete data, the alternating projection algorithms skip the projection onto the magnitude constraint in this region. Similarly, we only apply the relaxation $g_{\sigma} = \frac{1}{2\sigma} \sum_i \big||z_i| - b_i\big|^2$ on the known data for GPS. A simple exercise shows that
\begin{equation*}
    z^{k+1}_i = \begin{cases}
        \Big(\textrm{prox}_{t g_{\sigma}} \big( z^k - t \mathcal{F}y^k \big) \Big)_i & \mbox{if } b_i \textrm{ is known}, \\
        \big(z^k - t \mathcal{F}y^k\big)_i    & \mbox{otherwise.}
    \end{cases}
\end{equation*}

\begin{figure*}[h]
    \centering
    \includegraphics[height=3.6cm]{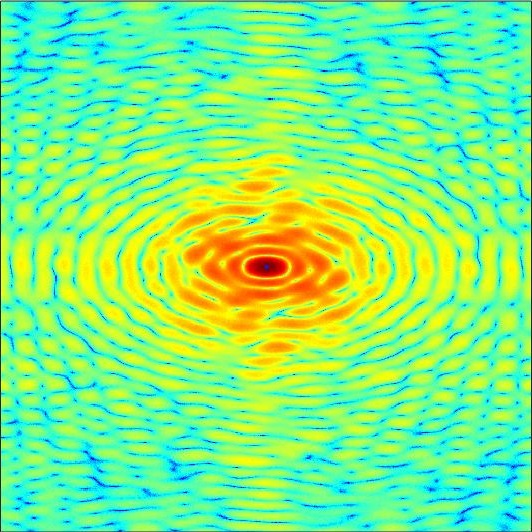} \hspace{1mm}
    \includegraphics[height=3.6cm]{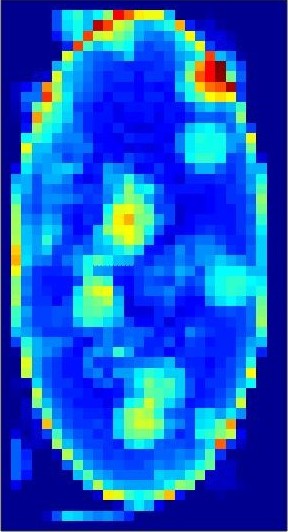} \hspace{1mm}
    \includegraphics[height=3.6cm]{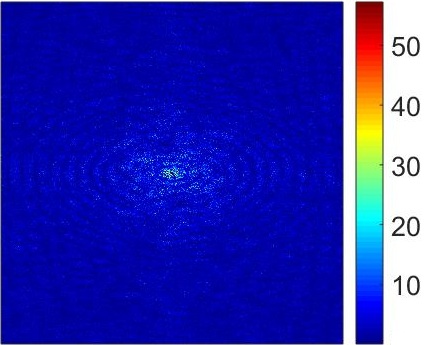} \\   
\vspace{2mm}
    \includegraphics[height=3.6cm]{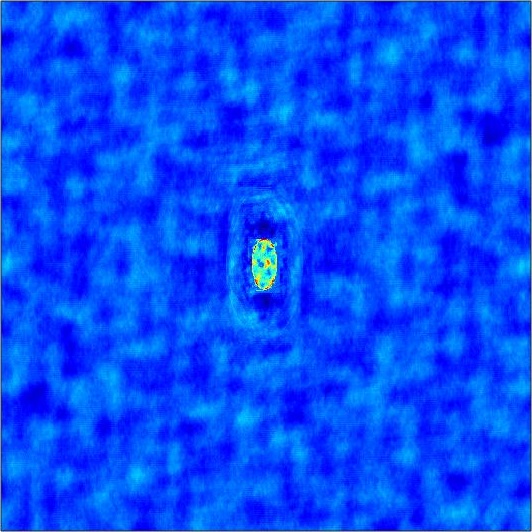} \hspace{1mm}
    \includegraphics[height=3.6cm]{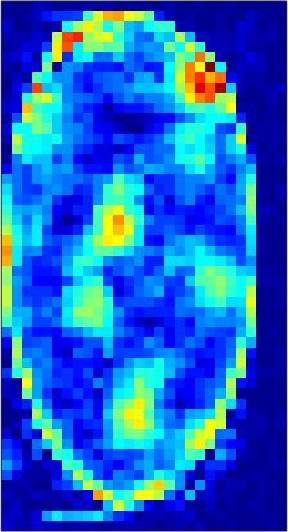} \hspace{1mm}
    \includegraphics[height=3.6cm]{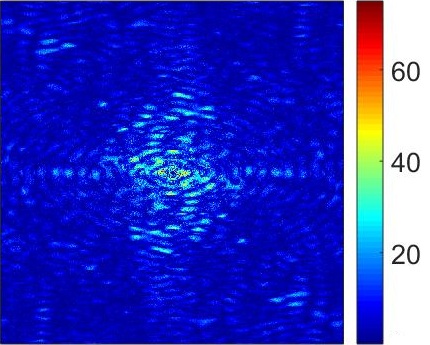}\\ \vspace{2mm}
    \includegraphics[height=3.6cm]{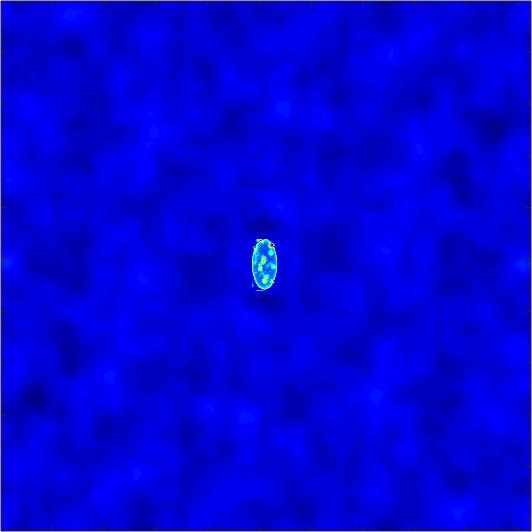} \hspace{1mm}
    \includegraphics[height=3.6cm]{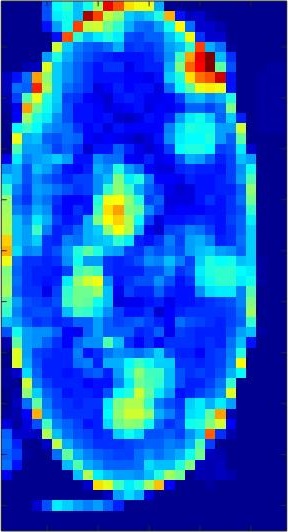} \hspace{1mm} 
    \includegraphics[height=3.6cm]{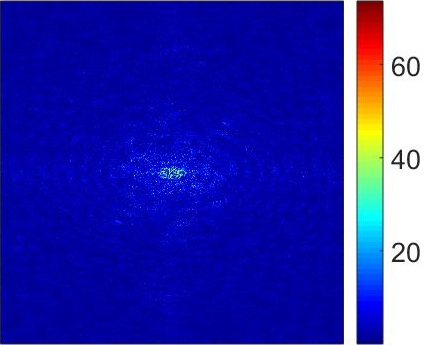}\\ \vspace{2mm}
    \includegraphics[height=3.6cm]{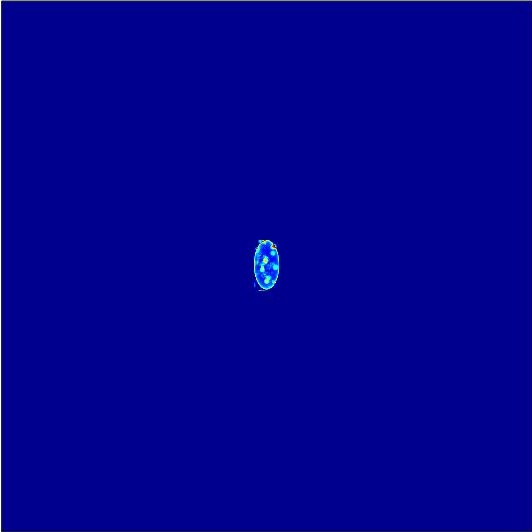} \hspace{1mm}
    \includegraphics[height=3.6cm]{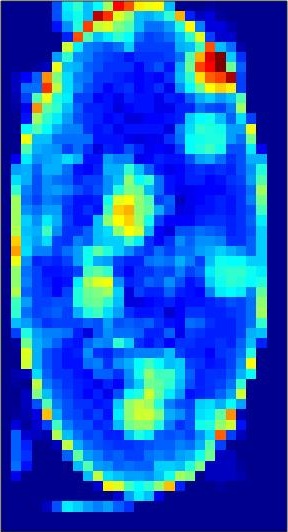} \hspace{1mm}
    \includegraphics[height=3.6cm]{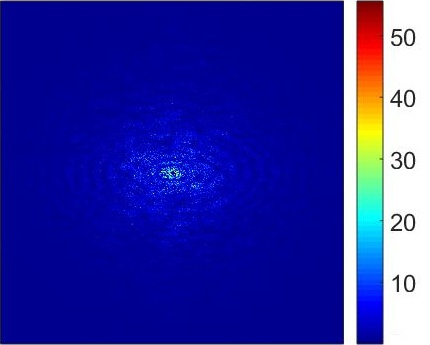} \\
    \vspace{2mm} 
    \includegraphics[height=3.6cm]{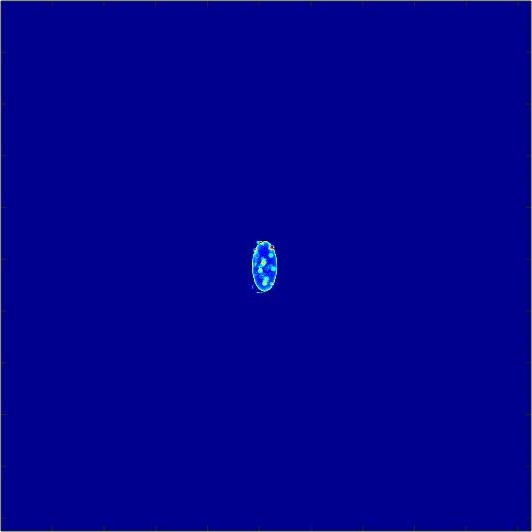} \hspace{1mm}
    \includegraphics[height=3.6cm]{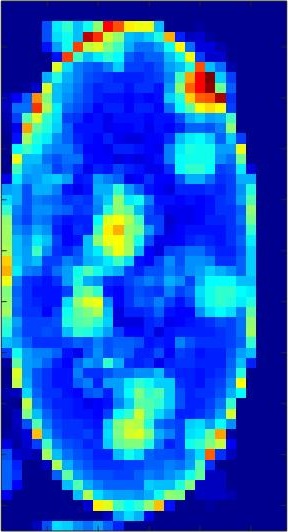} \hspace{1mm}
    \includegraphics[height=3.6cm]{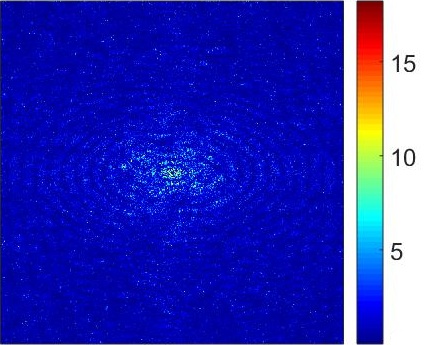}
    \caption{First row: vesicle model with log-scaled diffraction pattern (left) , zoom-in image (center) and residual (right). Second row: HIO: $R_F=12.87\%$, $R_{real}=21.14\%$. Third row: OSS: $R_F=6.08\%$,  $R_{real}=3.59\%$. Forth row: GPS-R $R_F=5.90\%$, $R_{real}=2.85\%$. Fifth row: GPS-F $R_F=5.89\%$ and $R_{real}=0.7\%$. Third column: residual.}
    \label{fig:vesicle}
\end{figure*}

%\bigskip

\section{Experimental Results}
\subsection{Reconstruction from simulated data}

% \begin{figure}[p!]
% 	\centering
%     flux = 1e7 \hspace{2.3cm} flux = 1e8 \hspace{2.3cm} flux = 1e9 \\
%     \vspace{1mm}
%     \includegraphics[width=4cm]{vesicle_DP_1e7.png} \hspace{1mm}
%     \includegraphics[width=4cm]{vesicle_DP_1e8.png} \hspace{1mm}
%     \includegraphics[width=4cm]{vesicle_DP_1e9.png}
%     \caption{simulated diffraction patterns of the vesicle model at three flux levels}
% \end{figure}

% \begin{table}[p!]
% \centering
% \begin{tabular}{|c|c|>{\columncolor[gray]{0.8}}c|c|>{\columncolor[gray]{0.8}}c|}
%     \hline
%     &\multicolumn{2}{c|}{$R_F$} & \multicolumn{2}{c|}{$R_{real}$} \\
%     \hline
%     & OSS & GPS-F & OSS & GPS-F \\
%     \hline
%     $\textrm{flux}=10^9$ & $1.673\%\pm 0.001\%$ & $1.426\%\pm 0.000\%$ & $1.745\%\pm 0.005\%$ & $0.135\%\pm 0.000\%$ \\
%     \hline
%     $\textrm{flux}=10^8$ & $4.709\%\pm 0.000\%$ & $4.634\%\pm 0.000\%$ & $2.185\% \pm 0.011\%$ & $1.770\%\pm 0.011\%$\\
%     \hline
%     $\textrm{flux}=10^7$ & $16.97\%\pm 0.001\%$ & $16.90\%\pm 0.004\%$ & $7.550\%\pm 0.021\%$ & $5.730\%\pm 0.008\%$\\
%     \hline
% \end{tabular}
% \caption{$R_F$ and $F_{real}$ as a function of flux (noise level) of vesicle data using OSS and GPS-F.}
% \label{table_flux}
% \end{table}
We expect GPS to be a reliable PR algorithm in the reconstruction of the images of weakly scattering objects, in particular biological specimens, which have become more popular \cite{Miao2003}. Since OSS has been shown to consistently outperform ER, HIO, ER-HIO, NR-HIO \cite{Rodriguez}, we perform both quantitative and qualitative comparisons between GPS and OSS. To simulate realistic experimental conditions, the Fourier magnitudes of a vesicle model are first corrupted with 5\% Poisson noise. $R_{noise}$ is defined to be the relative error with respect to the noise-free Fourier measurements
\begin{figure*}[h]
    \centering
    \includegraphics[height=5cm]{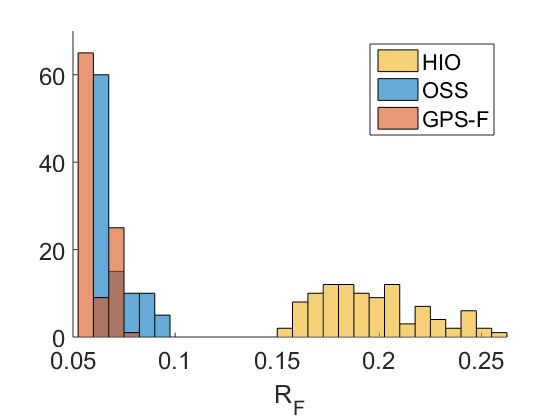} 
    \includegraphics[height=5cm]{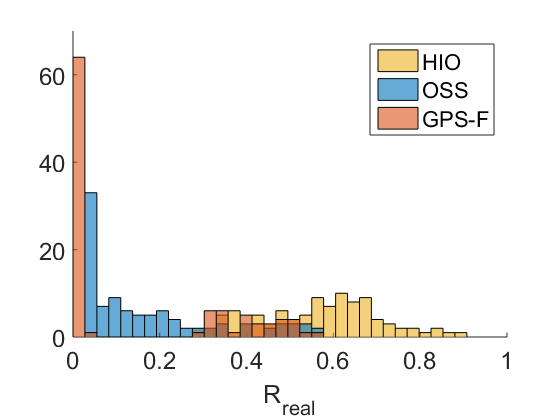}\\
    \includegraphics[height=5cm]{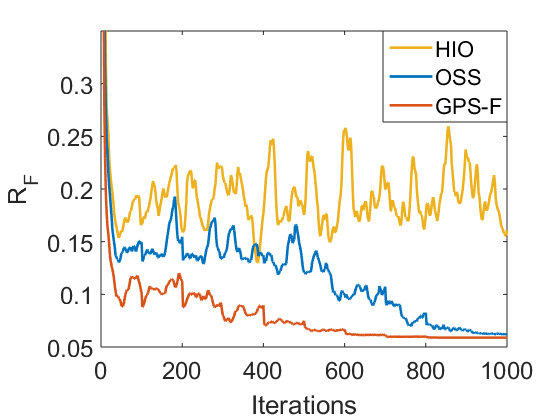}
    \includegraphics[height=5cm]{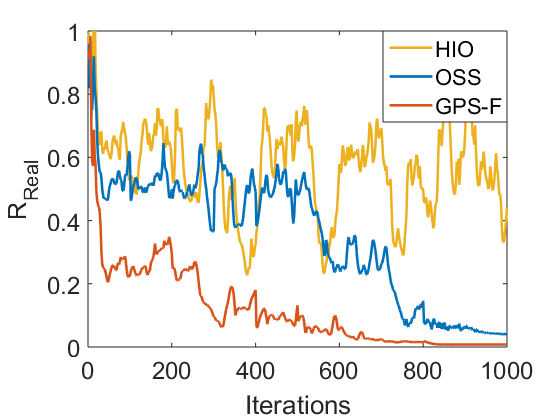}    
    \caption{Histogram (first row), and convergence (second row) of $R_F$ and $R_{real}$ on Vesicle data using HIO, OSS, GPS. GPS consistently produces smaller $R_F$ and $R_{real}$ than HIO or OSS. Moreover, GPS converges fastest with the fewest oscillations.}
    \label{vesicle_histogram}
\end{figure*}

\begin{equation*}
    R_{noise} = \frac{\sum_i \big| b_i - |\mathcal{F} u^o|_i \big|}{\sum_i |\mathcal{F} u^o |_i}
\end{equation*}
where $u^o$ is the noise-free model, and $b$ are noisy Fourier magnitudes. Due to the discrete nature of photon counting, experimentally measured data inherently contain Poisson noise that is directly related to the incident photon flux on the object. In addition to Poisson noise, the data is also corrupted with zero-mean Gaussian noise to simulate readout from a CCD. Any resulting negative values are set to zero. Therefore, an accurate simulation of $b_i$ can be calculated as 
\begin{equation*}
    b_i = \sqrt{ \textrm{Poisson} \bigg( \frac{ |\mathcal{F}u^o|_i^2 \cdot \textrm{flux}} {\|\mathcal{F}u^o\|^2} \bigg) + \mathcal{N}(0,\sigma) },
\end{equation*}
where $\sigma$ is proportional to the readout noise. \cite{Song:rw5052}
  
In some cases, the reconstructed image by an algorithm yields a small relative error $R_F$ but has low quality. This is the issue of over-fitting, an example of which can be seen in certain reconstructions using ER-HIO \cite{Rodriguez}. Smoothing is a technique to avoid data over-fitting. To validate results and show that our algorithm does not develop over-fitting, we measure the difference between the reconstructed image and the model by
\begin{equation*}
    R_{real} = \frac{\sum_i \big| u^k_i - u^o_i \big|}{\sum_i u^o_i}.
\end{equation*}

In addition, we also look at the residual $\mathrm{Res} = \|\mathcal{F}u| - b \|$ which measures the difference between the Fourier magnitudes of the reconstructed image and the experimental measurements. The residual can validate the noise-removal model, telling how much noise is removed. 
Figure \ref{fig:vesicle} shows the reconstruction of vesicle model from simulated noisy data using HIO, OSS, GPS-R, and GPS-F. GPS-F and GPS-R obtain lower $R_F$  and $R_{real}$ than OSS. Moreover, GPS-F can get very close to the ground truth with $R_F$ = $5.90\%$ and $R_{real}=0.7\%$. In addition to lower R values, GPS-R and GPS-F converge to zero outside the support. They both obtain lower residuals than OSS, specifically GPS-F produces the least residual. If we choose larger parameter for the $\ell_2$ gradient regularizer in Fourier space, we will get a smoother residual. Overall for realistic Gaussian and Poisson noise in Fourier space, GPS-F is a suitable noise-removal model.
%\sout{Since GPS-F produces less residual, for larger noise it may get trapped at local minima due to data over-fitting. In this case, we should use GPS-R and the combinations GPS-FR.}\\

Figure \ref{vesicle_histogram} shows the histogram and the convergence of $R_F$ and $R_{real}$ on 100 independent, randomly seeded runs using HIO, OSS and GPS on the simulated vesicle data. The histogram shows that GPS is more consistent and robust than OSS. It has a higher chance to converge to good minima with lower $R_F$ and $R_{real}$ than OSS. Furthermore, $R_F$ and $R_{real}$ of OSS scatter widely due to the initial value dependency. In contrast, GPS is more selective and less dependent on initial values. $R_F$ and $R_{real}$ of GPS are seen at a lower mean minimum with less variance.

Similar to HIO and ER, OSS keeps oscillating until a finer low-pass filter is applied. In contrast, GPS converges faster and is less oscillatory than OSS. In the presence of noise, alternating projection methods (ER, HIO, OSS) keep oscillating but do not converge. Applying smoothing and replacing the measurement constraint by the least squares fidelity term $g_{\sigma}(z) = \frac{1}{2\sigma} \||z|-b\|^2$ helps to reduce the oscillations; hence, the method can converge to a stationary point. Note that larger $\sigma$ reduces more oscillations, but also decreases the chance to escape from local minima. Alternating projection methods have $\sigma = 0$ since they impose measurement constraints. GPS obtains both smaller $R_F$, $R_{real}$, and lower variance. Even though $R_F$ are close to each other, $R_{real}$ of GPS is much smaller than OSS. This means GPS recovers the vesicle cell with higher quality than OSS. This simulation shows that GPS is more reliable and consistent than OSS. 

\subsection{Reconstructions from experimental data} 
\subsubsection{S. \textit{Pombe} Yeast Spore}
\begin{comment}
\begin{figure}[p!]
	\begin{subfigure}[b1]{0.2\textwidth}
    \includegraphics[height=4cm]{yeast_diff_log.jpg}
    \end{subfigure}
	\begin{subfigure}[b2]{\textwidth}
    \centering
    \hspace{-1cm} HIO \hspace{4cm} OSS \hspace{4cm} GPS \\
    \includegraphics[height=4.5cm]{yeast_HIO_mean.jpg} \hspace{.5mm}
    \includegraphics[height=4.5cm]{yeast_OSS_mean.jpg} \hspace{.5mm}
    \includegraphics[height=4.5cm]{yeast_GPS_mean.jpg}\\
    %\hspace{.1mm}
    \includegraphics[height=4.5cm]{yeast_HIO_var.jpg} \hspace{.5mm}
    \includegraphics[height=4.5cm]{yeast_OSS_var2.jpg}
    \hspace{.5mm}
    \includegraphics[height=4.5cm]{yeast_GPS_var.jpg}
    \end{subfigure}
    \caption{S. \textit{pombe} yeast spore log-scaled diffraction pattern, size $500 \times 500$ and the reconstructions by: HIO ($R_F = 15.697\% \pm 0.526\%$) OSS ($R_F = 9.775\% \pm 0.202\%$),  GPS ($R_F = 8.672\% \pm 0.025\%$). The means (first row) and variances (second row) are computed from the best 5 of 500 independent runs.}
    \label{yeast_result}
\end{figure}
\end{comment}
\begin{figure}[h]
	\centering
    \hspace{-1cm}\includegraphics[height=3.6cm]{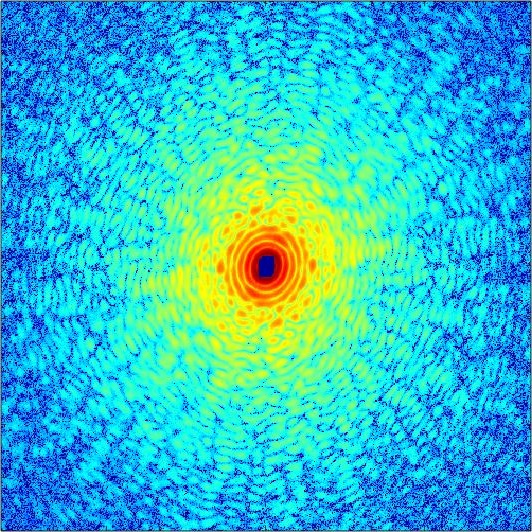} \\ \vspace{2mm}
    
    \includegraphics[height=3.6cm]{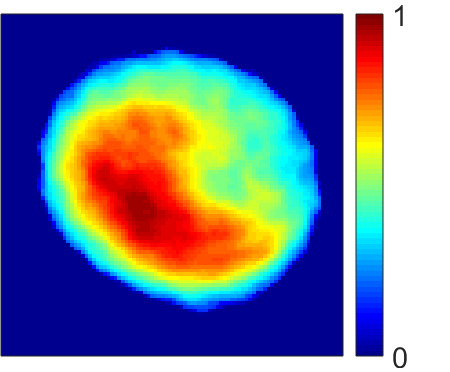} \hspace{-5mm}
    \includegraphics[height=3.6cm]{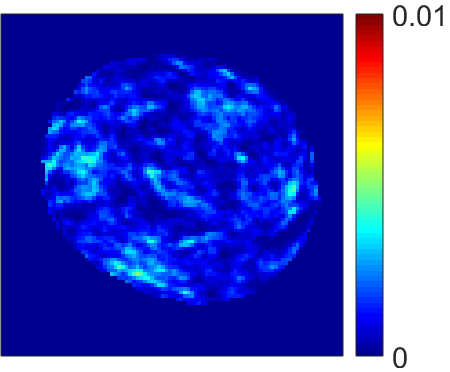} \\

    \includegraphics[height=3.6cm]{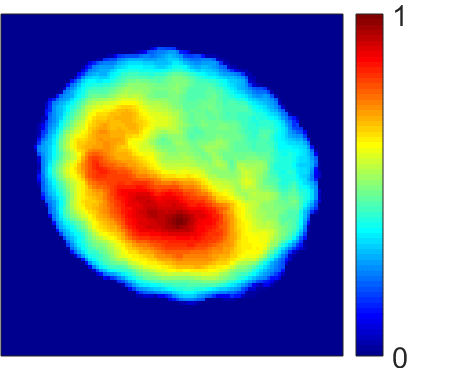} \hspace{-5mm}
    \includegraphics[height=3.6cm]{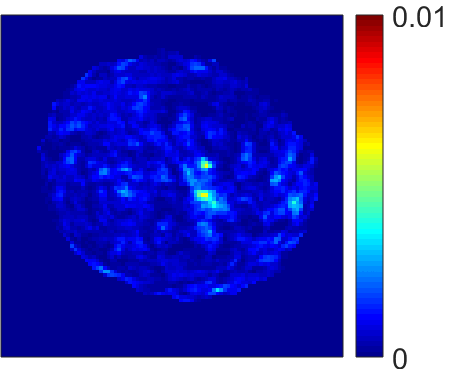}\\
    
\hspace{-1mm}
    \includegraphics[height=3.6cm]{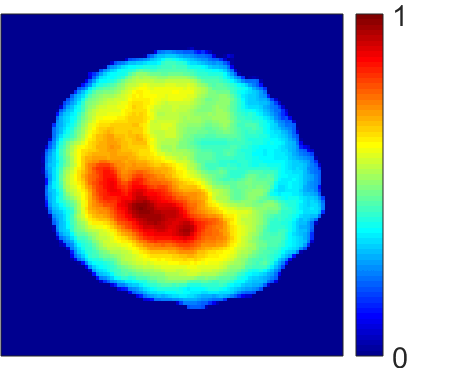} \hspace{-5mm}
    \includegraphics[height=3.6cm]{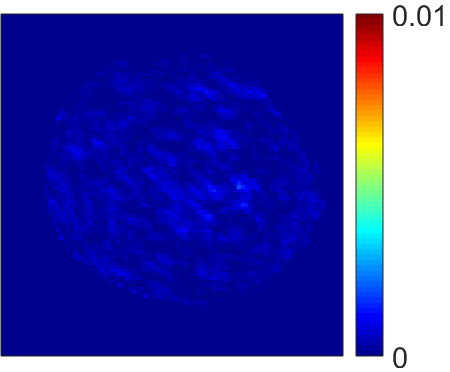}
    
    \caption{S. \textit{pombe} yeast spore log-scaled diffraction pattern, size $500 \times 500$ (top). Means (left column) and variance (right column) are computed from the best 5 of 500 independent reconstructions. HIO (first row): $R_F = 15.697\% \pm 0.526\%$, OSS (second row): $R_F = 9.775\% \pm 0.202\%$, and GPS (third row): $R_F = 8.672\% \pm 0.025\%$.}
    \label{yeast_result2}
\end{figure}

\begin{figure}[h]
    \centering
    \includegraphics[height=5cm]{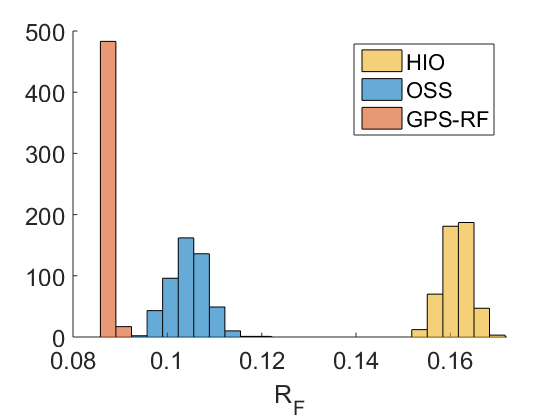}
    \includegraphics[height=5cm]{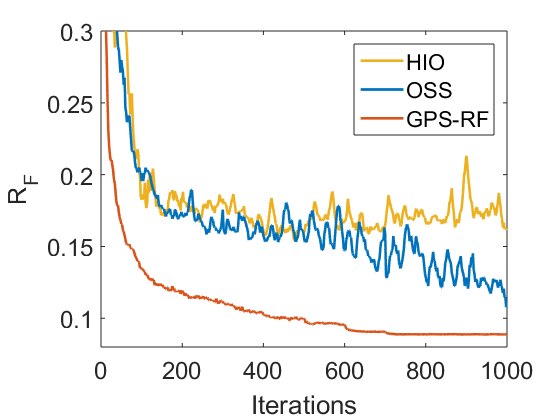}
    \caption{Top: histogram of $R_F$ in 500 independent runs (top). Bottom: the convergence curves of a singe construction of $R_F$ on S. \textit{pombe} yeast spore by GPS-RF, OSS, and HIO.}
    \label{yeast_hist_conv}
\end{figure}

To demonstrate the applicability of GPS to biological samples, we do phase retrieval on the diffraction pattern in figure \ref{yeast_result2} taken of a S. \textit{Pombe} yeast spore from an experiment done using beamline BL29 at SPring-8 \cite{Jiang11234}. We do 500 independent, randomly seeded reconstructions with each algorithm and record $R_F$, excluding the missing center. We choose default parameters for these experiments: $t=1$, $s=0.9$. The sequence of low-pass filters are chosen to be the same as in OSS \cite{Rodriguez}. For the first 400 iterations, $\sigma = 0.01$, then is increased to $\sigma = 0.1$ for the remaining 600 iterations. The left column of figure \ref{yeast_result2} is the mean of the best 5 reconstructions obtained by the respective algorithm. The right column shows the variance of the same 5 reconstructions. It is evident from the variance that GPS achieves more consistent reconstructions. Figure \ref{yeast_hist_conv} shows the histogram and convergence of $R_F$. We can conclude that not only are GPS-R results the most consistent, but also the most faithful to the experimental data.\\ 

\subsubsection{Nanorice}

\begin{comment}
\begin{figure}[p!]
    \centering
    %vesicle \hspace{2.5cm} 
    nanorice1 \hspace{2cm} nanorice2 \hspace{2.3cm} yeast spore \\ \vspace{1mm}
    %\includegraphics[height=4cm]{vesicle_diffpat.jpg}
    \includegraphics[height=4cm]{nanorice1_diffpat.jpg} \hspace{1mm}
    \includegraphics[height=4cm]{nanorice2_diffpat.jpg} \hspace{1mm}
    \includegraphics[height=4cm]{yeast_diff_log.jpg}
    \caption{Diffraction pattern of nanorice1 $253 \times 253$, nanorice2 $253 \times 253$, and log-scale S. \textit{pombe} yeast spore $500\times 500$ with missing data in the center.}
    \label{diffpat}
\end{figure}
\end{comment}

\begin{figure}[h]
    \centering
    nanorice1 \hspace{2cm} nanorice2 \\ \vspace{2mm}
    \includegraphics[height=3.6cm]{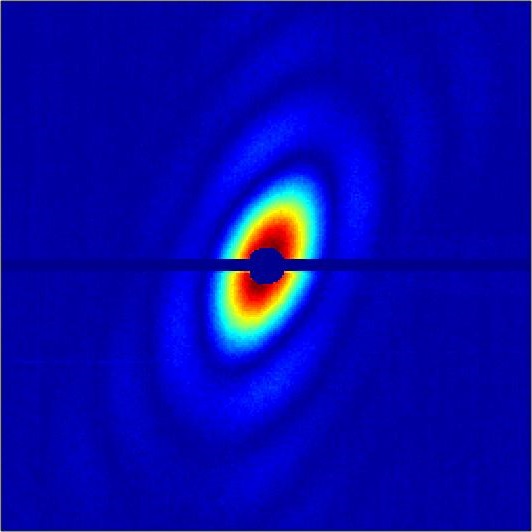} \hspace{1mm}
    \includegraphics[height=3.6cm]{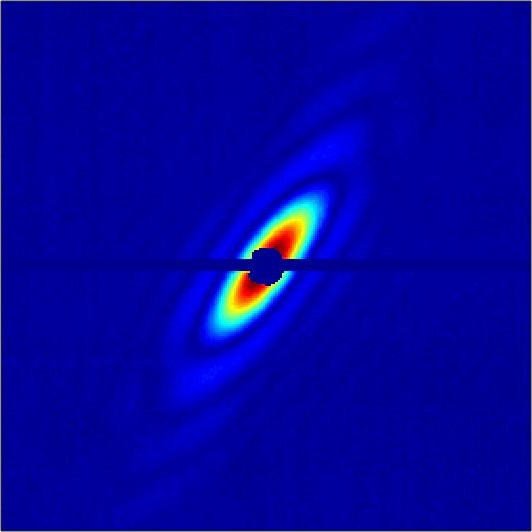}\\
    \vspace{2mm}
    \includegraphics[width=3.6cm]{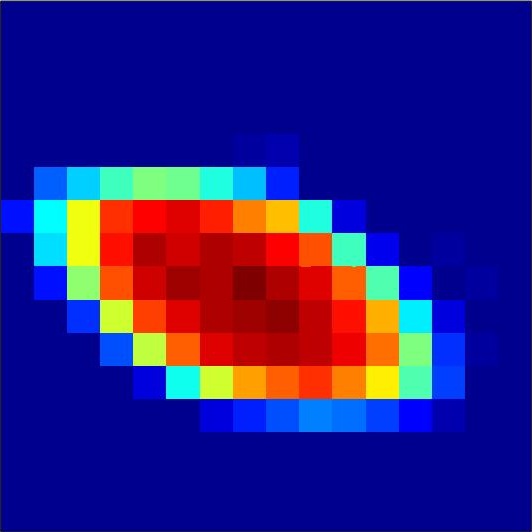}  \hspace{1mm}
    \includegraphics[width=3.6cm]{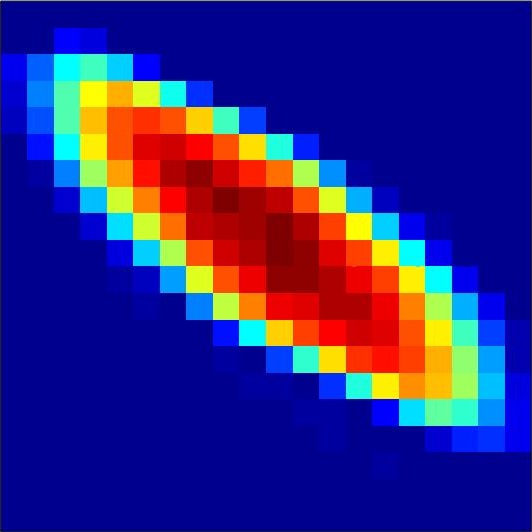} \\ \vspace{2mm} %\hspace{.1mm}
    \includegraphics[width=3.6cm]{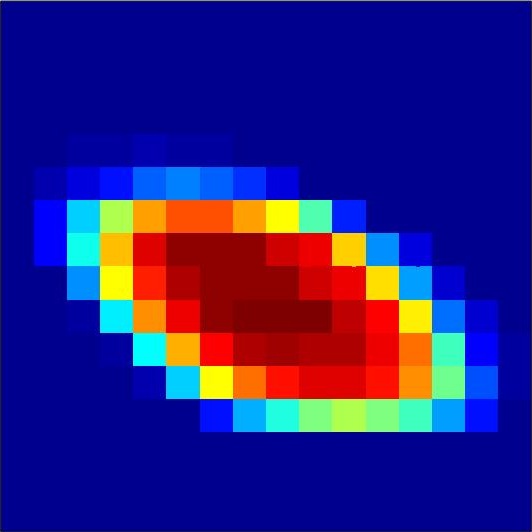}  \hspace{1mm}
    \includegraphics[width=3.6cm]{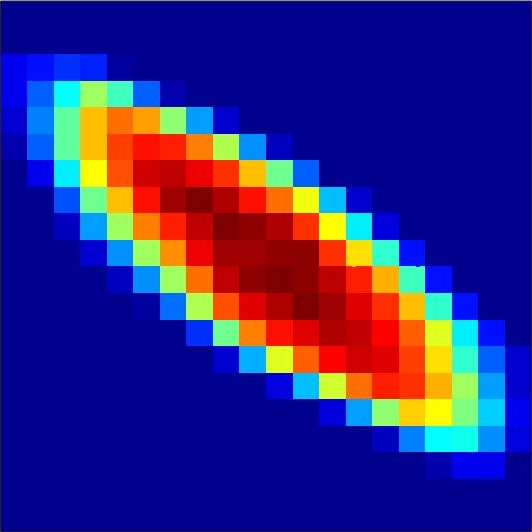} 
     \caption{Diffraction pattern of nanorice1 and nanorice2 $253 \times 253$ (first row) and reconstructions using OSS (second row: $R_F=18.23\%,\, 16.32\%$) and GPS-F (third row: $R_F=17.40\%,\, 15.83\%$) respectively. GPS obtains less noise on the boundary and lower relative error $R_F$.}
    \label{nanorice}
\end{figure}

% $nanorice1: OSS (R_F=18.23\%$), GPS-F ($\sigma=1$) ($R_F=17.42\%$). nanorice2 OSS ($R_F=16.32\%$), GPS-F ($\sigma=1$) ($R_F=15.86\%$)

\begin{figure*}[h]
    \centering
    \includegraphics[height=5cm]{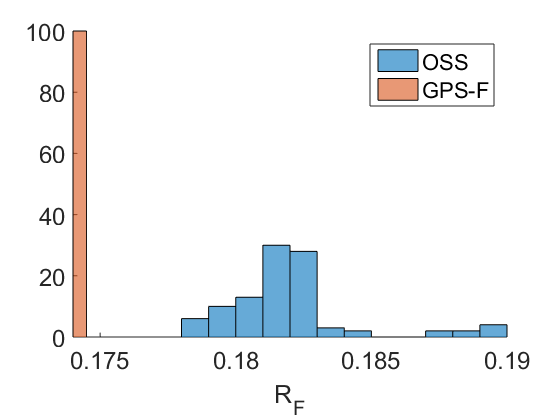}
    \includegraphics[height=5cm]{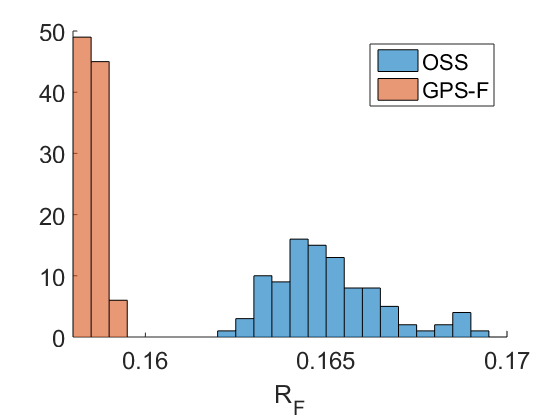}\\
    \includegraphics[height=5cm]{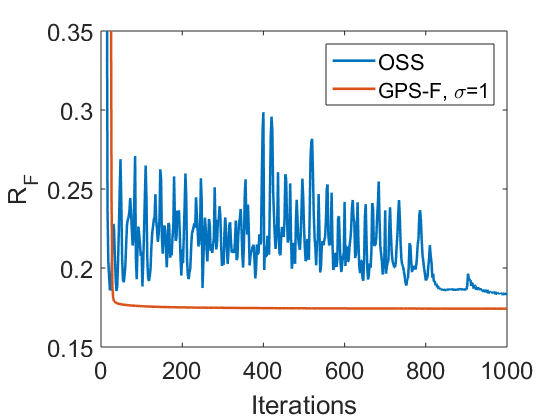}
    \includegraphics[height=5cm]{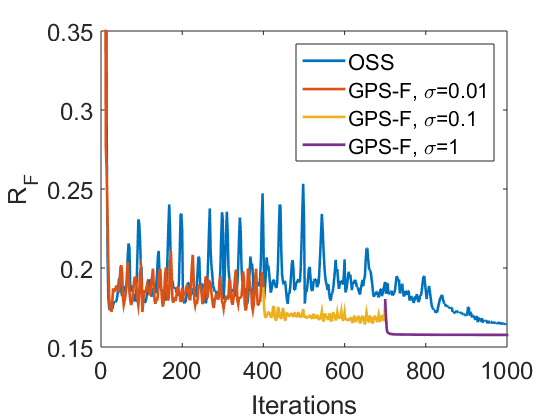}
    \caption{Histogram (first column) and convergence (second column) of OSS and GPS-F on nanorice1 (first row) and nanorice2 (second row). The results of HIO are omitted due to lack of convergence.}
    \label{nanorice_hist_conv}
\end{figure*}

To demonstrate the generality of GPS, we also do testing with experimental data of inorganic samples. The diffraction patterns shown in the top row of figure \ref{nanorice} from ellipsoidal iron oxide nanoparticles (referred to as \lq{nanorice1}\rq and \lq{nanorice2}\rq) were taken at the AMO instrument at LCLS at an X-ray energy of 1.2keV \cite{Kassemeyer:12}. This data is freely available online on the CXIDB \cite{maia_coherent_2012}. We choose default parameters for these experiments: $t=1$, $s=0.9$. The sequence of low-pass filters are chosen to be the same as in OSS \cite{Rodriguez}. The fidelity parameter $\sigma$ is chosen small for the first 800 iterations, specifically $\sigma \in [0,0.01]$, to produce oscillations which is necessary for the algorithm to skip bad local minima. Once the reconstruction arrives at a good local minimum region, we increase $\sigma$ to reduce oscillations. This later value of $\sigma$ depends on noise level and data. We test different values of $\sigma$ and $\sigma = 1$ has been found to be the optimal for both nanorice data. Figure \ref{nanorice} shows OSS(second row) and GPS-F(third row) reconstructions of the two nanorice particles. Figure \ref{nanorice_hist_conv} shows again that GPS obtains more consistent and faithful reconstructions as compared to those obtained by OSS. GPS-F with $\sigma = 1$ converges to lower relative error than OSS at all times. OSS cannot get lower relative error because $\sigma = 0$ does not work for this case. In general, alternating projection methods do not treat noise correctly. For example, in this case, HIO keeps oscillating but does not converge. Therefore, its results are omitted here. The better approach, OSS model, can reduce oscillations by smoothing but this is not enough. In contrast, the least squares $g_{\sigma}(z) = \frac{1}{2\sigma} \| |z| - b \|^2$ of GPS works for noise removal since relaxing the constraints allows GPS to reach lower relative error. The values of $\sigma$ depend on noise level and type. To optimize the convergence of GPS-F on nanorice2, we apply $\sigma=0.01$ for the first 400 iterations, $\sigma=0.1$ for the next 300 iterations, and $\sigma=1$ for the last 300 iterations. This test shows the effect of $\sigma$ on the convergence. OSS ($\sigma=0$) oscillates the most. GPS with $\sigma=0.01, 0.1, 1$ oscillates less and less. As $\sigma$ increases, GPS also gets to lower $R_F$. The algorithm finally reaches a stable minimum as $\sigma$ goes up to 1. Continuing to increase $\sigma$ does not help with $R_F$. Choosing large $\sigma$ in the beginning may reduce oscillations but also limit the mobility to skip local minima.  We recommend start with small $\sigma$ and then gradually increase it until the iterative process reaches a stable minimum. \\
 
\section{Conclusion}
In this work, we have developed a fast and robust phase retrieval algorithm GPS for the reconstruction of images from noisy diffraction intensities. Similar to \cite{binaryrelax}, the Moreau-Yosida regularization was used to relax the hard constraints considered in the noiseless model.
GPS utilizes a primal-dual algorithm and a noise-removal technique, in which the $\ell_2$ gradient smoothing is effectively applied on either real or Fourier space of the dual variable. GPS shows more reliable and consistent results than OSS, HIO for the reconstruction of weakly scattered objects such as biological specimens. Looking forward, we aim to explore the role of dual variables in non-convex optimization. Smoothing the dual variable, which is equivalent to smoothing the gradient of convex conjugate, represents a new and effective technique that can in principle be applied to other non-smooth, non-convex problems. \\

%\noindent \textbf{Acknowledgement.} 

\begin{comment}
\begin{figure}
    \begin{subfigure}{.33\textwidth}
        \centering
        flux=1e7
        \includegraphics[width=4cm]{vesicle_DP_1e7.png}
    \end{subfigure}\hfill
    \begin{subfigure}{.33\textwidth}
        \centering
        flux=1e8
        \includegraphics[width=4cm]{vesicle_DP_1e8.png}
    \end{subfigure}
    \begin{subfigure}{.33\textwidth}
        \centering
        flux=1e9
        \includegraphics[width=4cm]{vesicle_DP_1e9.png}
    \end{subfigure}
    \caption{simulated diffraction patterns of the vesicle model at three flux levels}
    \label{vesicle_DP_sim}
\end{figure}
\end{comment}

\addtolength{\textheight}{-12cm}   % This command serves to balance the column lengths

%%%%%%%%%%%%%%%%%%%%%%%%%%%%%%%%%%%%%%%%%%%%%%%%%%%%%%%%%%%%%%%%%%%%%%%%%%%%%%%%

%%%%%%%%%%%%%%%%%%%%%%%%%%%%%%%%%%%%%%%%%%%%%%%%%%%%%%%%%%%%%%%%%%%%%%%%%%%%%%%%

%%%%%%%%%%%%%%%%%%%%%%%%%%%%%%%%%%%%%%%%%%%%%%%%%%%%%%%%%%%%%%%%%%%%%%%%%%%%%%%%
%\section*{APPENDIX}


\section*{ACKNOWLEDGMENT}

We thank Professor Jose A. Rodriguez for help and providing experimental data of the yeast spore.  This work was supported by STROBE: A National Science Foundation Science $\&$ Technology Center, under Grant No. DMR 1548924 as well as DOE-DE-SC0013838. 

%%%%%%%%%%%%%%%%%%%%%%%%%%%%%%%%%%%%%%%%%%%%%%%%%%%%%%%%%%%%%%%%%%%%%%%%%%%%%%%%
\clearpage
\bibliographystyle{unsrt}
\bibliography{main}

\end{document}